\documentclass[12pt]{amsart}

\usepackage{amssymb}
\usepackage{amsmath}
\usepackage{amscd}
\usepackage{float}
\usepackage{graphicx}
\usepackage{amsfonts}
\usepackage{pb-diagram}

\newtheorem{thm}{Theorem}
\newtheorem{cor}{Corollary}
\newtheorem{lem}{Lemma}
\newtheorem{prop}{Proposition}
\newtheorem{defn}{Definition}
\newtheorem{rem}{Remark}

\begin{document}

\title{An invariant for singular knots}

\author{J. Juyumaya}
\address{Departamento de Matem\'aticas, Universidad de Valpara\'{\i}so \\
Gran Breta\~na 1091, Valpara\'{\i}so, Chile.}
\email{juyumaya@uvach.cl}

\author{S. Lambropoulou}
\address{ Departament of Mathematics,
National Technical University of Athens,
Zografou campus, GR-157 80 Athens, Greece.}
\email{sofia@math.ntua.gr}
\urladdr{http://www.math.ntua.gr/$\tilde{~}$sofia}

\thanks{Both authors were partially supported by Fondecyt 1085002, NTUA and  Dipuv.}



\date{}
\maketitle
\begin{abstract}
In this paper we introduce a Jones-type invariant for singular knots, using a Markov trace on the Yokonuma--Hecke algebras ${\rm Y}_{d,n}(u)$ and the theory of singular braids. The Yokonuma--Hecke algebras have a natural topological interpretation in the context of framed knots. Yet, we show that there is a homomorphism of the singular braid monoid $SB_n$ into the algebra  ${\rm Y}_{d,n}(u)$. Surprisingly, the trace does not normalize directly to yield a singular link invariant, so a condition must be imposed on the trace variables. Assuming this condition, the invariant satisfies a skein relation involving singular crossings, which arises from a quadratic relation in the algebra ${\rm Y}_{d,n}(u)$.
\end{abstract}

\section{Introduction}


A {\it singular link on $n$ components} is the image of a smooth immersion of $n$ copies of the circle in $S^3$, that has finitely many singularities, called {\it singular crossings}, which are all ordinary double points.  So, a singular link is like a classical link, but with a finite number of transversal self-intersections permitted. A singular link on one component is a {\it singular knot}. Some examples of singular knots and links are given in Figure~\ref{fig5}. We shall say `knots' throughout meaning `knots and links'.

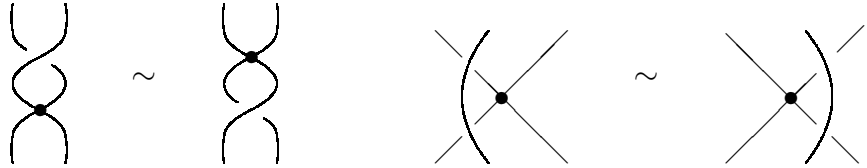
\begin{figure}[H]
\begin{center}
 

\begin{picture}(320,48)

\qbezier(10,50)(9,55)(10,60)
\qbezier(30,50)(31,55)(30,60)
\qbezier(10,0)(9,5)(10,10)
\qbezier(30,0)(31,5)(30,10)

\qbezier(10,10)(10,15)(20,20)
\qbezier(20,20)(30,25)(30,30)
\qbezier(30,10)(30,15)(20,20)
\qbezier(20,20)(10,25)(10,30)

\qbezier(10,30)(10,35)(20,40)
\qbezier(20,40)(30,45)(30,50)
\qbezier(30,30)(30,34)(25,37)
\qbezier(15,43)(10,46)(10,50)


\put(55,30){$\sim$}

\qbezier(90,50)(89,55)(90,60)
\qbezier(110,50)(111,55)(110,60)
\qbezier(90,0)(89,5)(90,10)
\qbezier(110,0)(111,5)(110,10)

\qbezier(90,10)(90,15)(100,20)
\qbezier(100,20)(110,25)(110,30)
\qbezier(110,10)(110,14)(105,17)
\qbezier(95,23)(90,26)(90,30)
\qbezier(90,30)(90,35)(100,40)
\qbezier(100,40)(110,45)(110,50)
\qbezier(110,30)(110,35)(100,40)
\qbezier(100,40)(90,45)(90,50)
\put(17.5,17){$\bullet$}

\put(97.5,37){$\bullet$}
\put(170,0){\line(1,1){10}}
\put(185,15){\line(1,1){35}}

\put(220,0){\line(-1,1){35}}
\put(180,40){\line(-1,1){10}}

\qbezier(190,0)(170,25)(190,50)

\put(245,30){$\sim$}
\put(280,0){\line(1,1){34}}
\put(322,42){\line(1,1){10}}

\put(330,0){\line(-1,1){10}}
\put(314,15){\line(-1,1){34}}
\qbezier(310,0)(330,25)(310,50)
\put(192,21.5){$\bullet$}
\put(301.5,21.5){$\bullet$}

\end{picture}
\caption{The singular Reidemeister moves}\label{fig0}
\end{center}
\end{figure}

Two singular links $K_1, K_2$ are {\it isotopic}, that is, topologically equivalent, if there is an orientation preserving self-homeomorphism of $S^3$ carrying one to the other, such that it preserves a small rigid disc around each singular crossing (rigid-vertex isotopy). In terms of diagrams, $K_1, K_2$ are isotopic if and only if any two diagrams of theirs differ by planar isotopy and a finite sequence of the classical and the singular Reidemeister moves. In Figure~\ref{fig0} we illustrate the main two singular Reidemeister moves. The others are the obvious variants of these, with different crossings.

\smallbreak

By their definition, singular links may admit a well-defined orientation on each component. Then the isotopy moves are considered with all possible orientations. It is a well--known fact that every isotopy invariant $\mathcal L$ of classical oriented links extends to an invariant of singular oriented links, by means of the rule:
\[
{\mathcal L}(L_\times) = {\mathcal L}(L_+) - {\mathcal L}(L_-)
\]
where $L_+, L_-$ and $L_\times$ are identical diagrams, except for the place of one crossing, where it is positive, negative or singular respectively (see Figure \ref{fig4}).

\smallbreak

A {\it singular braid on $n$ strands} is the image of a smooth immersion of $n$ arcs  in $S^3$, that has finitely many singularities, the singular crossings, which are all ordinary double points, such that the ends are arranged into $n$ collinear top endpoints and into $n$ collinear bottom endpoints and such that there are no local maxima or minima. So, a  singular braid is like a classical braid, but with a finite number of singular crossings allowed.

\begin{figure}[H]
\begin{picture}(200,70)
\put(-22,67){\tiny{$1$}}
\qbezier(-19,15)(-19,45)(-19,65)
\put(-16,38){$\ldots$}
\qbezier(0,30)(0,35)(10,40)
\qbezier(10,40)(20,45)(20,50)
\qbezier(20,30)(20,34)(15,37)
\qbezier(5,43)(0,46)(0,50)

\qbezier(0,50)(0,70)(0,60)
\qbezier(20,50)(20,70)(20,60)
\qbezier(0,30)(0,10)(0,20)
\qbezier(20,30)(20,10)(20,20)

\put(38,67){\tiny{$n$}}
\qbezier(40,15)(40,45)(40,65)
\put(21,38){$\ldots$}
\put(160,67){\tiny{$1$}}
\qbezier(161,15)(161,45)(161,65)
\put(165,38){$\ldots$}


\qbezier(200,30)(200,35)(190,40)
\qbezier(190,40)(180,45)(180,50)
\qbezier(180,30)(180,35)(190,40)
\qbezier(190,40)(200,45)(200,50)

\qbezier(200,50)(200,70)(200,60)
\qbezier(180,50)(180,70)(180,60)
\qbezier(200,30)(200,10)(200,20)
\qbezier(180,30)(180,10)(180,20)

\put(219,67){\tiny{$n$}}
\qbezier(221,15)(221,45)(221,65)
\put(202,38){$\ldots$}
\put(5, 5){$\sigma_{i}$}
\put(-1,67){\tiny{$i$}}
\put(13,67){\tiny{$i+1$}}
\put(187.5, 36.5){$\bullet$}
\put(185, 5){$\tau_{i}$}
\put(179,67){\tiny{$i$}}
\put(192,67){\tiny{$i+1$}}
\end{picture}
\caption{The elementary braids $\sigma_i$ and $\tau_i$}
\label{fig1}
\end{figure}
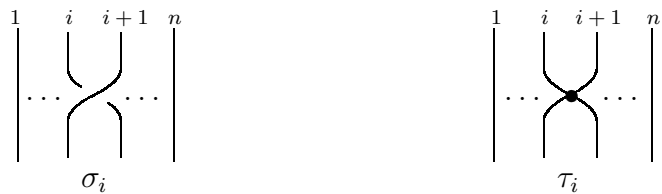

 Two singular braids are isotopic if there is a rigid-vertex isotopy taking one to the other, which fixes the endpoints of the strands and preserves the braid structure.
Algebraically, the set of singular braids on $n$ strands, denoted $SB_n$, forms a monoid with the usual concatenation of braids, the so-called {\it singular braid monoid}. It was introduced in different contexts by Baez\cite{ba}, Birman\cite{bi} and Smolin\cite{sm}.
 $SB_n$ is generated by the unit, by the classical elementary braids $\sigma_i$ with their inverses, and by the corresponding elementary singular braids $\tau_i$ (view Figure~\ref{fig1}):
$$
1, \sigma_1, \ldots,
\sigma_{n-1}, \sigma_1^{-1}, \ldots , \sigma_{n-1}^{-1},\tau_1, \ldots ,  \tau_{n-1}
$$
which satisfy the relations below. These reflect precisely the singular braid isotopy.
\begin{equation}\label{relsbn}
\begin{array}{rclcll}
\sigma_i\sigma_i^{-1}    & = & \sigma_i^{-1} \sigma_i   & = & 1 &  \mbox{for all}\,\, i \\

[\sigma_i, \sigma_j] & = & [\sigma_i, \tau_j] & = & [\tau_i, \tau_j] =0 &
\mbox{for}\, \vert i-j\vert  >1 \\

[\sigma_i, \tau_i] & = & 0 & & & \mbox{for all}\,\, i \\
\sigma_i\sigma_j\sigma_i & = & \sigma_j\sigma_i\sigma_j &   &   &
 \mbox{for}\,  \vert i-j\vert = 1\\
 \sigma_i \sigma_j \tau_i &  = &   \tau_j \sigma_i \sigma_j
& & &{\rm for}\,  \vert i-j\vert = 1.
\end{array}
\end{equation}

\smallbreak

The singular braid monoid  $SB_n$ embeds in a group, the singular braid group, see \cite{fkr}.


\smallbreak

The {\it closure} of a singular braid is defined like the ordinary closure of a classical braid, whereby we join the endpoints of corresponding strands by simple arcs. The closure of a singular braid $\omega$  shall be denoted $\widehat{\omega}$.
 In analogy to the classical setting, oriented singular links may be isotoped to closed singular braids. For a proof of Alexander's theorem for singular links see \cite{bi}.

 Let now $\cup_n SB_n$ denote the inductive limit associated to the natural monomorphisms of monoids $SB_n \hookrightarrow SB_{n+1}$. In analogy to the Markov theorem for classical braids, Gemein  proved in \cite{ge1} the following result (compare also with \cite{la} for an $L$-move version).

\begin{thm}[Gemein, 1997]\label{gema}
 Two singular braids in $\cup_n SB_n$ have isotopic closures if and only if they differ by singular braid relations and a finite sequence of the following moves:

\vspace{.07in}
\noindent (i) \, Real conjugation: \ $\sigma_i \omega \sim \omega \sigma_i, \quad \omega,\sigma_i \in SB_n$

\vspace{.07in}
\noindent (ii) \ Singular commuting: \ $\tau_i\omega \sim \omega \tau_i, \quad \omega,\tau_i \in SB_n$

\vspace{.07in}
\noindent (iii) Real stabilization: \ $\omega \sim \omega \sigma_n^{\pm 1}, \quad \omega \in SB_n$\
\end{thm}

\begin{figure}[H]
\begin{center}
 
\begin{picture}(210,95)


\qbezier(0,30)(0,50)(0,70)
\qbezier(70,30)(70,50)(70,70)
\qbezier(0,30)(35,30)(70,30)
\qbezier(0,70)(35,70)(70,70)

\qbezier(140,30)(140,50)(140,70)
\qbezier(210,30)(210,50)(210,70)
\qbezier(140,30)(175,30)(210,30)
\qbezier(140,70)(175,70)(210,70)


\qbezier(25,5)(25,10)(35,15)
\qbezier(35,15)(45,20)(45,25)
\qbezier(45,5)(45,10)(35,15)
\qbezier(35,15)(25,20)(25,25)

\qbezier(25,0)(25,2.5)(25,5)
\qbezier(45,0)(45,2.5)(45,5)
\qbezier(25,25)(25,27.5)(25,30)
\qbezier(45,25)(45,27.5)(45,30)

\put(32,11){$\bullet$}

\qbezier(165,75)(165,80)(175,85)
\qbezier(175,85)(185,90)(185,95)
\qbezier(185,75)(185,80)(175,85)
\qbezier(175,85)(165,90)(165,95)

\qbezier(165,95)(165,97.5)(165,100)
\qbezier(185,95)(185,97.5)(185,100)
\qbezier(165,70)(165,72.5)(165,75)
\qbezier(185,70)(185,72.5)(185,75)

\put(173,82){$\bullet$}

\put(30,45){$\omega$}
\put(170,45){$\omega$}
\put(100,45){$\sim$}


\qbezier(12.5,0)(12.5,15)(12.5,30)
\qbezier(57.5,0)(57.5,15)(57.5,30)

\qbezier(12.5,70)(12.5,85)(12.5,100)
\qbezier(57.5,70)(57.5,85)(57.5,100)

\qbezier(27,70)(27,85)(27,100)
\qbezier(43,70)(43,85)(43,100)

\qbezier(152.5,0)(152.5,15)(152.5,30)
\qbezier(197.5,0)(197.5,15)(197.5,30)

\qbezier(167,0)(167,15)(167,30)
\qbezier(183,0)(183,15)(183,30)

\qbezier(152.5,70)(152.5,85)(152.5,100)
\qbezier(197.5,70)(197.5,85)(197.5,100)
\end{picture}
\caption{Singular commuting}\label{fig2}
\end{center}
\end{figure}
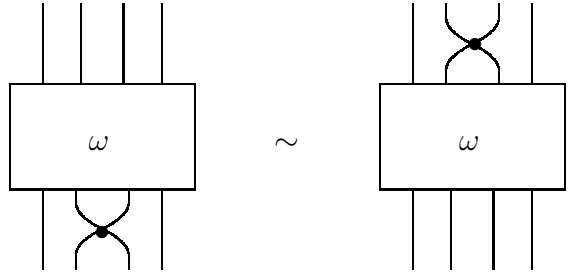

Moves (i) and (iii) of Theorem~\ref{gema} are the two well--known moves of the classical Markov theorem for classical braids. Move (ii) is illustrated in Figure~\ref{fig2}.

\smallbreak

Using the Alexander and Markov theorems for singular links and braids, it is possible to construct singular link invariants via  Markov traces on quotient algebras. Baez\cite{ba} introduced the Vassiliev algebra as the quotient of $\Bbb C \,SB_n \otimes \Bbb C(\epsilon)$ by the ideal generated by the expressions $\sigma_i - \sigma_i^{-1}-\epsilon\, \tau_i$, which give rise to the following relations in the algebra:
$$
g_i - g_i^{-1} = \epsilon\, \tau_i
$$
He then showed that $\Bbb C$-valued Vassiliev--Gussarov invariants are in one-to-one correspondence with homogeneous (in $\epsilon$) Markov traces on the algebra.

\smallbreak
 More recently, Paris and Rabenda\cite{para} defined the singular Hecke algebra as the quotient of the algebra $\Bbb C(q)[SB_n]$ by the ideal generated by the  well--known quadratic relations $\sigma_i^2 -(q-1) \sigma_i - q\, 1$  of the classical Iwahori--Hecke algebra of type $A$, which give rise to the following relations in the algebra:
$$
g_i - q\, g_i^{-1} = (q-1)\, 1
$$
 They also constructed on these algebras singular Markov traces, from which --upon normalization, according to Theorem~1-- they derived the universal HOMFLYPT (2-variable Jones polynomial) analogue for oriented singular links. Kauffman and Vogel\cite{kavo} constructed analogues of the HOMFLYPT and the Kauffman polynomial for singular links by using diagrammatic methods. For their HOMFLYPT analogue it is shown in \cite{para} that it is a specialization of the universal HOMFLYPT  for singular links (but, as observed in \cite{para}, this specialization does not make much difference).

\bigbreak

In the present paper we construct an invariant for singular links using the Yokonuma--Hecke algebras and Markov traces defined on them. The  Yokonuma--Hecke algebra ${\rm Y}_{d,n}(u)$ can be defined as a quotient of the modular framed braid group algebra ${\Bbb C} {\mathcal F}_{d,n}$ (classical framed braids with framings modulo $d$) by the quadratic relations $g_i^2 = 1 + (u-1)e_i(1 - g_i)$ (here also $\sigma_i$ corresponds to $g_i$). The elements $e_i$ are certain idempotents in ${\Bbb C} {\mathcal F}_{d,n}$ and they are expressions of the framing generators $t_i, t_{i+1}$, see (\ref{ei}). Further, in \cite{ju} Juyumaya constructed a linear Markov trace on the  Yokonuma--Hecke algebra ${\rm Y}_{d,n}(u)$, that is, a linear trace which supports the Markov property: ${\rm tr}(ag_n) = z\,{\rm tr}(a)$ for $z \in {\Bbb C}$ and for any $a \in {\rm Y}_{d,n}(u)$, see Theorem~\ref{trace}. For details and topological interpretations of the above we refer the reader to \cite{jula} and \cite{jula2}, where the algebras ${\rm Y}_{d,n}(u)$ and the traces in \cite{ju} are used in the context of classical and $p$-adic framed braids and framed links.

\smallbreak
Here, we first map homomorphically the singular braid monoid $SB_n$ into the algebra ${\rm Y}_{d,n}(u)$ via the  map:
$$
\begin{array}{cccl}
\delta : & SB_n & \longrightarrow & {\rm Y}_{d,n}(u) \\
 & \sigma_i  & \mapsto & g_i \\
 & \tau_i    & \mapsto & p_i = e_i(1-g_i)
\end{array}
$$
In the image $\delta (SB_n)$ the following relations hold:
\begin{equation}\label{gipi}
g_i-g_i^{-1} = (u^{-1}-1)p_i
\end{equation}
 For an illustration of the corresponding quadratic relations (Eq. \ref{cure}) see Figure~\ref{fig3}.
 We then consider  the Markov trace on ${\rm Y}_{d,n}(u)$ constructed
in \cite{ju}. So, we obtain a map from the singular braid monoid $SB_n$ to the complex numbers.
 A normalization of this trace according to the singular braid equivalence of Theorem~\ref{gema} should yield a singular link invariant. But this turns out not to be the case (not even for the restriction to classical braids and links). In fact, we need to impose a condition on the variables of the trace, the `$E$--condition', see  Definition \ref{defcon}. Surprisingly, there are non-trivial sets of complex numbers satisfying the $E$--condition. Given now the $E$--condition we normalize the trace to obtain an invariant $\Delta$ of singular links (Theorem \ref{invariant}).

\smallbreak
To the best of our knowledge, the algebra ${\rm Y}_{d,n}(u)$ is the first example of an algebra, which can admit so different topological interpretations: in the context of framed braids as well as in the context of singular braids. Also, the  Markov trace in \cite{ju} is the only Markov trace we know of that does not normalize directly to yield a link invariant. Finally, it is worth adding that, for modulus $d$ equal to $1$, we have $e_i=1$ and the algebra ${\rm Y}_{1,n}(u)$ coincides with the classical Iwahori--Hecke algebra of type $A$. Also, the trace in  \cite{ju} coincides with the Ocneanu trace \cite{jo}.

\smallbreak
We would like to thank the referee of the paper for very useful comments and remarks.

\section{The Yokonuma--Hecke algebra and a Markov trace}

\subsection{\it Relations in ${\rm Y}_{d,n}(u)$}

We fix a $u \in {\Bbb C}\backslash \{0,1\}$.  The Yokonuma--Hecke algebra, denoted by ${\rm Y}_{d,n}(u)$, is a ${\Bbb C}$--associative algebra  generated by the elements
$$
1, g_1, \ldots, g_{n-1}, t_1, \ldots, t_{n}
$$
subject to the following relations:
$$
\begin{array}{rcll}
g_ig_j & = & g_jg_i & \mbox{for $\vert i-j\vert > 1$}\\
g_ig_jg_i & = & g_jg_ig_j & \mbox{for $ \vert i-j\vert = 1$}\\
t_i t_j & =  &  t_j t_i &  \mbox{for all $ i,j$}\\
t_j g_i & = & g_i t_{s_i(j)} & \mbox{for all $ i,j$}\\
t_j^d   & =  &  1 & \mbox{for all $j$}
\end{array}
$$
where $s_i(j)$ is the result of applying the transposition $s_i=(i, i+1)$ to $j$, together with the extra quadratic relations:
\begin{equation}\label{quadr}
g_i^2 = 1 + (u-1) \, e_{i} - (u-1) \, e_{i} \, g_i \qquad \mbox{for all $i$}
\end{equation}
 where
\begin{equation}\label{ei}
e_i :=\frac{1}{d}\sum_{m=0}^{d-1}t_i^m t_{i+1}^{-m}
\end{equation}
The first four relations are defining relations for the classical framed braid group, with the $t_j$'s being interpreted as the `elementary framings' (framing 1 on the $j$th strand). The relations $t_j^d = 1$ mean that the framing of each strand is regarded modulo~$d$. So, the algebra ${\rm Y}_{d,n}(u)$ arises naturally  as a quotient of the modular framed braid group algebra over the quadratic relations (\ref{quadr}). But in the present paper we shall give a different topological interpretation to ${\rm Y}_{d,n}(u)$, in relation to singular knots and links.

It is easily verified that the elements $e_i$ are idempotents.  Also, that the elements $g_i$ are invertible in ${\rm Y}_{d,n}(u)$. Indeed:
\begin{equation}\label{invrs}
g_i^{-1} = g_i - (u^{-1} - 1)\, e_i + (u^{-1} - 1)\, e_i\, g_i
\end{equation}

As noted in the Introduction, $d=1$ implies $e_i=1$ and $p_i=1-g_i$. So $g_i^2 = u+ (1-u) g_i$, and the algebra ${\rm Y}_{1,n}(u)$ coincides with the Iwahori--Hecke algebra of type $A$.
For more details on the algebra ${\rm Y}_{d,n}(u)$ and for further topological interpretations see \cite{jula,jula2} and references therein. In ${\rm Y}_{d,n}(u)$ we have the following relations.

\begin{lem}\label{leme}
For the elements $e_i$ and for $1\leq i,j\leq n-1$ the following relations hold:
$$
\begin{array}{rcll}
e_ie_j & = & e_j e_i \\
e_ig_i & = & g_i e_i \\
e_i g_j & = & g_je_i &  {\text for } \vert i - j \vert >1 \\
e_j g_ig_j & = & g_ig_je_i & {\text for } \vert i - j \vert =1
\end{array}
$$
\end{lem}

\begin{proof}
The first three claims are easy to check (see Lemma 4 and Proposition 5 in \cite{jula}). We will check the last one. Let $j= i+1$. From the defining relations in the algebra ${\rm Y}_{d,n}(u)$ we have:

\noindent $t_{i+1}g_ig_{i+1}= g_it_{i}g_{i+1} = g_ig_{i+1}t_{i}$. Similarly, $t_{i+2}g_ig_{i+1} = g_ig_{i+1}t_{i+1}$. Then
$$
e_{i+1}g_ig_{i+1} = \frac{1}{d}\sum_{m=0}^{d-1} t_{i+1}^mt_{i+2}^{-m}g_ig_{i+1} =\frac{1}{d}\sum_{m =0}^{d-1} g_ig_{i+1}t_{i}^{m}t_{i+1}^{-m} =g_ig_{i+1}e_{i}.
$$
 The proof for $j= i-1$ is completely analogous.
\end{proof}

\subsection{\it The elements $p_i$}

For $1 \leq i\leq  n-1$ we define the elements $p_i \in {\rm Y}_{d,n}(u)$ by the formula:
\begin{equation}\label{piai}
p_i = e_i(1-g_i)
\end{equation}
Then, the quadratic relations (\ref{quadr}) in the algebra ${\rm Y}_{d,n}(u)$ may be rewritten as:
\begin{equation}\label{cure}
g_i^2 = 1 + (u-1)p_i
\end{equation}

\begin{prop}\label{relp}
For the elements $p_i$ and for $1\leq i,j\leq n-1$, we have the relations:
$$
\begin{array}{rcll}
e_ip_i & = & p_ie_i = p_i \\
p_i^k & = & (u+1)^{k-1}p_i & {\text for } k\in {\Bbb N} \\
g_ip_i & = & p_ig_i = -up_i &  \\
g_ip_j & = & p_jg_i & {\text for }\vert i - j \vert >1 \\
p_ip_j & =  & p_jp_i & {\text for } \vert i - j \vert >1 \\
p_j g_ig_j & =  & g_ig_jp_i & {\text for } \vert i - j \vert =1
\end{array}
$$
\end{prop}

\begin{proof}
The proofs follow from Eq.~\ref{piai}, from  Lemma~\ref{leme} and by direct computations. For example, we shall check the second relation. For $k=2$ we have $p_i^2 = e_i(1-g_i)e_i(1-g_i) = e_i^2(1-g_i)^ 2 = e_i(1-2g_i + g_i^2)= e_i\left(1-2g_i + 1 + (u-1)p_i\right) = 2e_i (1-g_i) + (u-1)e_ip_i$. Then, by the first relation we have $p_i^ 2  = 2p_i + (u-1)p_i= (u+1)p_i$. For any $k>2$ we apply induction.
\end{proof}

Note that the elements $p_i$ are not invertible in ${\rm Y}_{d,n}(u)$. We now define for fixed $a \in {\Bbb C}$ the following map.

\begin{equation}\label{deltaa}
\begin{array}{cccl}
\delta_{a} : & SB_n & \longrightarrow & {\rm Y}_{d,n}(u) \\
 & \sigma_i  & \mapsto & ag_i \\
 & \tau_i    & \mapsto & p_i
\end{array}
\end{equation}
In particular, we shall denote:
$$
\delta := \delta_1
$$

\begin{thm}\label{homom}
The map $\delta_{a}$ defines a monoid homomorphism.
\end{thm}

\begin{proof}
The proof follows immediately by comparing relations (\ref{relsbn}) in $SB_n$ with the relations in Proposition~\ref{relp}.
\end{proof}

\subsection{\it Topological interpretations}

\begin{figure}[H]
\begin{center}
 
\begin{picture}(200,70)

\qbezier(0,60)(-1,65)(0,70)
\qbezier(20,60)(21,65)(20,70)
\qbezier(0,10)(-1,15)(0,20)
\qbezier(20,10)(21,15)(20,20)

\qbezier(0,20)(0,25)(10,30)
\qbezier(10,30)(20,35)(20,40)
\qbezier(20,20)(20,24)(15,27)
\qbezier(5,33)(0,36)(0,40)
\qbezier(0,40)(0,45)(10,50)
\qbezier(10,50)(20,55)(20,60)
\qbezier(20,40)(20,44)(15,47)
\qbezier(5,53)(0,56)(0,60)
\put(40,40){$=$}

\put(65,10){\line(0,1){62}}
\put(85,10){\line(0,1){62}}
\put(103,40){$+ \quad (u-1)$}

\qbezier(200,30)(200,35)(190,40)
\qbezier(190,40)(180,45)(180,50)
\qbezier(180,30)(180,35)(190,40)
\qbezier(190,40)(200,45)(200,50)
\put(187.5,37){$\bullet$}

\qbezier(200,50)(202,60)(202,70)
\qbezier(180,50)(178,60)(178,70)
\qbezier(202,10)(202,20)(200,30)
\qbezier(178,10)(178,20)(180,30)

\end{picture}
\caption{The quadratic relations involve singular crossings}\label{fig3}
\end{center}
\end{figure}
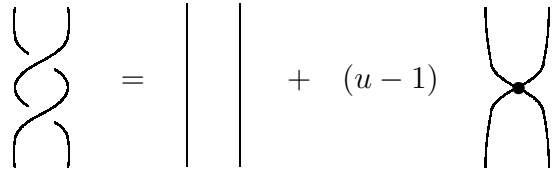

 We shall now give topological interpretations for the elements of the subalgebra $\delta (SB_n)$ of the algebra ${\rm Y}_{d,n}(u)$. By Theorem~\ref{homom}, monomials in $g_i, g_i^{-1}, p_i$ may be viewed as singular braids, such that $g_i, g_i^{-1}$ correspond respectively to $\sigma_i, \sigma_i^{-1}$  (for $a=1$) and $p_i$ corresponds to the singular crossing $\tau_i$. These elementary braids are subject to the quadratic relations in Eq.~\ref{cure}. These relations are illustrated in Figure \ref{fig3}, where, for simplicity, we omit the identity strands.  Multiplying Eq. \ref{cure} by $g_i^{-1}$ and using Proposition~\ref{relp} we obtain the  equivalent relation (\ref{gipi}).

\smallbreak

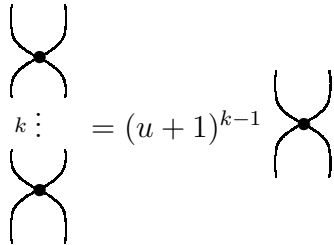
\begin{figure}[H]
\begin{center}
\begin{picture}(120,85)

\qbezier(10,0)(9,5)(10,10)
\qbezier(30,0)(31,5)(30,10)
\qbezier(10,10)(10,15)(20,20)
\qbezier(20,20)(30,25)(30,30)
\qbezier(30,10)(30,15)(20,20)
\qbezier(20,20)(10,25)(10,30)
\qbezier(10,30)(9,32.5)(10,35)
\qbezier(30,30)(31,32.5)(30,35)
\qbezier(10,55)(9,57.5)(10,60)
\qbezier(30,55)(31,57.5)(30,60)
\qbezier(10,60)(10,65)(20,70)
\qbezier(20,70)(30,75)(30,80)
\qbezier(30,60)(30,65)(20,70)
\qbezier(20,70)(10,75)(10,80)
\qbezier(10,80)(9,85)(10,90)
\qbezier(30,80)(31,85)(30,90)


\put(40,40){$={\tiny (u+1)^{k-1}}$}

\qbezier(110,55)(109,60)(110,65)
\qbezier(130,55)(131,60)(130,65)
\qbezier(110,25)(109,30)(110,35)
\qbezier(130,25)(131,30)(130,35)

\qbezier(110,35)(110,40)(120,45)
\qbezier(120,45)(130,50)(130,55)
\qbezier(130,35)(130,40)(120,45)
\qbezier(120,45)(110,50)(110,55)
\put(17.5,17){$\bullet$}
\put(17.5,67){$\bullet$}
\put(117.5,42){$\bullet$}
\put(11, 42){\tiny{$k$}}
\put(18.5, 40){$\vdots$}

\end{picture}
\caption{The relation $p_i^k = (u+1)^{k-1}p_i$}\label{fig6}
\end{center}
\end{figure}

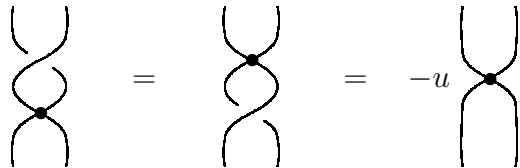
\begin{figure}[H]
\begin{center}
 

\begin{picture}(200,60)

\qbezier(10,50)(9,55)(10,60)
\qbezier(30,50)(31,55)(30,60)
\qbezier(10,0)(9,5)(10,10)
\qbezier(30,0)(31,5)(30,10)

\qbezier(10,10)(10,15)(20,20)
\qbezier(20,20)(30,25)(30,30)
\qbezier(30,10)(30,15)(20,20)
\qbezier(20,20)(10,25)(10,30)

\qbezier(10,30)(10,35)(20,40)
\qbezier(20,40)(30,45)(30,50)
\qbezier(30,30)(30,34)(25,37)
\qbezier(15,43)(10,46)(10,50)
\put(55,30){$=$}

\qbezier(90,50)(89,55)(90,60)
\qbezier(110,50)(111,55)(110,60)
\qbezier(90,0)(89,5)(90,10)
\qbezier(110,0)(111,5)(110,10)

\qbezier(90,10)(90,15)(100,20)
\qbezier(100,20)(110,25)(110,30)
\qbezier(110,10)(110,14)(105,17)
\qbezier(95,23)(90,26)(90,30)
\qbezier(90,30)(90,35)(100,40)
\qbezier(100,40)(110,45)(110,50)
\qbezier(110,30)(110,35)(100,40)
\qbezier(100,40)(90,45)(90,50)
\put(135,30){$=\quad-u$}
\qbezier(180,23)(180,28)(190,33)
\qbezier(190,33)(200,38)(200,43)
\qbezier(200,23)(200,28)(190,33)
\qbezier(190,33)(180,38)(180,43)
\qbezier(180,43)(179,55)(180,60)
\qbezier(200,43)(201,55)(200,60)
\qbezier(180,0)(179,5)(180,23)
\qbezier(200,0)(201,5)(200,23)
\put(17.5,17){$\bullet$}
\put(97.5,37){$\bullet$}

\put(187.5,30){$\bullet$}

\end{picture}
\caption{The relation $g_ip_i = p_ig_i = -up_i$}\label{fig7}
\end{center}
\end{figure}

\begin{figure}[H]
\begin{center}
 
\begin{picture}(230,80)
\put(50,80){\vector(-1,-1){50}}
\put(28,52){\vector(1,-1){22}}
\put(0,80){\line(1,-1){22}}

\put(118,58){\line(1,1){22}}
\put(112,52){\vector(-1,-1){22}}
\put(90,80){\vector(1,-1){50}}

\put(230,80){\vector(-1,-1){50}}
\put(180,80){\vector(1,-1){50}}

\put(202,52){$\bullet$}
\put(15,10){$L_+$}
\put(110,10){$L_{-}$}
\put(200,10){$L_{\times}$}

\end{picture}
\caption{$L_{+}$, $L_{-}$ and $L_{\times}$}\label{fig4}
\end{center}
\end{figure}
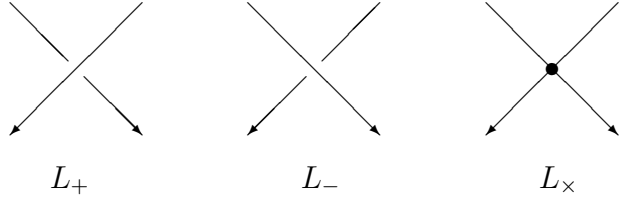

Beyond relations (\ref{relsbn}), the images of the generators of $SB_n$ under the map $\delta$ are also subject to the extra relations $p_i^k = (u+1)^{k-1}p_i$ and $p_ig_i = -up_i$ of Proposition~1. These are illustrated in Figures~\ref{fig6} and \ref{fig7}, where the identity strands are also omitted.

\smallbreak

\smallbreak

Note that in this topological set--up there are no obvious interpretations for the generators $t_i$ and the elements  $e_i$.

\subsection{\it A Markov trace on ${\rm Y}_{d,n}(u)$}

Let now $\cup_{n}{\rm Y}_{d,n}(u)$ denote the inductive limit associated to the natural inclusions ${\rm Y}_{d,n}(u) \subset {\rm Y}_{d,n+1}(u)$. In \cite{ju} the following theorem is proved.

\begin{thm}[Juyumaya, 2004]\label{trace}
Let $z$, $x_1$, $\ldots$, $x_{d-1}$ be in ${\Bbb C}$.
There exists a unique linear map ${\rm tr}$ on $\cup_{n}{\rm Y}_{d,n}(u)$ with
values  in  ${\Bbb C}$ satisfying the rules:
$$
\begin{array}{rcll}
{\rm tr}(ab) & = & {\rm tr}(ba) & \\
{\rm tr}(1) &  = & 1 & \\
{\rm tr}(ag_n) & = & z\,{\rm tr}(a)& \qquad (a \in {\rm Y}_{d,n}(u))\\
{\rm tr}(at_{n+1}^m) & = & x_m{\rm  tr}(a) & \qquad ( a \in {\rm Y}_{d,n}(u), \, 1\leq m\leq d-1).\\
\end{array}
$$
\end{thm}

As noted in the Introduction, for $d=1$ the trace restricts to the first three rules and it coincides with Ocneanu's trace on the Iwahori--Hecke algebra, which was used to construct the 2--variable Jones polynomial for classical knots and links, see \cite{jo}.

\section{The $E$--condition and an invariant for singular knots}

In view of Theorems~\ref{gema}, \ref{homom} and \ref{trace} we would like to construct an isotopy invariant for singular knots and links. According to Theorem~\ref{gema}, such an invariant has to agree on the singular links $\widehat{\omega}$, $\widehat{\omega\sigma_n}$ and $\widehat{\omega\sigma_n^{- 1}}$, for any $\omega \in SB_n$. Now, having present the recipe of Jones\cite{jo} for constructing the (2-variable) Jones polynomial for classical knots, we will try to define an invariant by re-scaling and normalizing the trace $\rm tr$. By Eq.~\ref{invrs} we have:
$$
{\rm tr}(\omega g_n^{-1}) = {\rm tr}(\omega g_n) - (u^{-1}-1){\rm tr}(\omega e_n)
 + (u^{-1}-1){\rm tr}(\omega e_ng_n)
$$
In order that the invariant agrees on the closures of the braids $\omega{\sigma_n}^{-1}$ and $\omega{\sigma_n}$ we need that ${\rm tr}(\omega g_n^{-1})$ factorizes through ${\rm tr} (\omega)$.  For the first term we have: ${\rm tr}(\omega g_n) = z\, {\rm tr}(\omega)$. Further:
\begin{equation}\label{egomega}
{\rm tr} (\omega e_ng_n)=\frac{1}{d}\sum_{m=0}^{d-1}{\rm tr}(\omega t_n^mt_{n+1}^{-m} g_n )=\frac{1}{d}\sum_{m=0}^{d-1}z\, {\rm tr}(\omega) = z\,{\rm tr}(\omega)
\end{equation}
since ${\rm tr}(\omega t_n^mt_{n+1}^{-m} g_n)= {\rm tr}( \omega t_n^mg_n t_{n}^{-m})= z\,{\rm tr}( \omega t_n^mt_n^{-m}) = z\,{\rm tr}(\omega)$.

\subsection{\it The $E$--condition}

 From the above analysis it is clear that $\rm tr$ needs  to satisfy also the following multiplicative property:
\begin{equation}\label{mul}
{\rm tr}(\omega e_n) = {\rm tr}( e_n)\, {\rm tr}(\omega)
\end{equation}
Unfortunately, we do not have such a  nice formula for ${\rm tr}(\omega \, e_n)$. The underlying reason on the framed braid level (that is, for the natural interpretation for elements in ${\rm Y}_{d,n}(u)$) is that $e_n$ involves the $n$th strand of $\omega$.
Yet, by imposing some conditions on the indeterminates $x_i$ it is possible to have property (\ref{mul}). Before giving these conditions let us define the following elements in $Y_{d,n}(u)$:
$$
e_i^{(m)} : =\frac{1}{d}\sum_{s=0}^{d-1}t_i^{m + s} t_{i+1}^{-s}  \quad  \text{and} \quad
e_i := e_i^{(0)}
$$
Also, the corresponding elements in ${\Bbb C}[z, x_1, \ldots, x_{d-1}]$:
$$
\zeta^{(m)} : = \frac{1}{d}\sum_{s=0}^{d-1}x_{s+ m}x_{d-s} = {\rm tr}(e_i^{(m)}) \quad  \text{and} \quad
\zeta := \zeta^{(0)} = {\rm tr}(e_i)
$$
where the sub-indices of the indeterminates  are regarded modulo $d$.

\begin{defn}\label{defcon}\rm
 We shall say that the set ${\mathcal X}_d :=\{x_1,\ldots , x_{d-1}\}$ of complex numbers satisfies  the {\it $E$--condition} if it satisfies the following system of $d-1$ non--linear equations  in ${\Bbb C}$:
$$
\zeta^{(m)} = x_m \zeta \qquad (1\leq m \leq d-1)
$$
Or, equivalently:
\begin{equation}\label{Esystem}
\sum_{s=0}^{d-1}x_{m + s} x_{d-s}  =  x_m \sum_{s=0}^{d-1} x_{s} x_{d-s} \qquad (1\leq m \leq d-1)
\end{equation}
where the sub-indices on the $x_j$'s are regarded modulo $d$ and $x_0 =x_d:=1$.
\end{defn}

Surprisingly, there exist non--trivial sets ${\mathcal X}_d$ satisfying the $E$--condition. For example, taking $x_i= \theta^i$, where $\theta$ is a primitive $d$th root of  unity. We note that for this solution we have $\zeta ={\rm tr}(e_j) =1$ and ${\rm tr}(p_j) =1 - z$. For $d=3,4$ and $5$ we run the Mathematica program and we found  other  solutions of the $E$--system, for which:
$$
{\rm tr}(e_j) \neq 1.
$$
 For example in the case $d=3$, where we have the $E$--system:
$$
\begin{array}{lcr}
x_1 + x_2^2& = & 2 x_1^2 x_2 \\
x_1^2 + x_2& = & 2 x_1 x_2^2
\end{array}
$$
we have the non-trivial solutions:
$$
x_1=x_2 = -\frac{1}{2}\quad \text{or}\quad
x_1=\frac{1}{3}\left(\frac{1}{3}- \frac{3i\sqrt{3}}{4}\right), \
x_2=\frac{1}{4}\left( 1 + i\sqrt{3}\right)
$$
Also, the solution where we take the conjugates in the previous one.
Another more interesting example is the set formed by the elements
$$
x_i :=\frac{-(-1)^{i(d-1)}}{d-1}\qquad (1\leq i \leq d-1)
$$
We then have $\zeta = {\rm tr}(e_j) = 1/(d-1)$. For explanations about the somewhat `mysterious' $E$--condition and for a thorough discussion on the solutions of the $E$--system we refer the reader to \cite{jula2}.

\subsection{\it A singular link invariant}

We are now close to our aim. Indeed, assuming the $E$--condition we have the following.

\begin{thm}\label{thmcon}
If ${\mathcal X}_d$ satisfies the $E$--condition, then for all $\omega \in
{\rm Y}_{d,n}(u)$ we have
$$
{\rm tr}(\omega e_{n}) = {\rm tr}( e_{n}) \, {\rm tr}(\omega) = \zeta \, {\rm tr}(\omega).
$$
\end{thm}
\begin{proof}
See \cite{jula2}.
\end{proof}

\begin{cor}\label{zetapn}
If ${\mathcal X}_d$ satisfies the $E$--condition, then for all $\omega \in
{\rm Y}_{d,n}(u)$ we have
$$
{\rm tr}(\omega p_{n}) = {\rm tr}( p_{n})\, {\rm tr}(\omega) = (\zeta-z) \, {\rm tr}(\omega).
$$
\end{cor}

\begin{proof}
By (\ref{mul}) we have: ${\rm tr}(\omega p_{n}) = {\rm tr}(\omega e_n(1-g_n)) = {\rm tr}(\omega e_n) - {\rm tr}(\omega e_ng_n)$. So, by Theorem~\ref{thmcon} and by (\ref{egomega}): ${\rm tr}(\omega p_{n}) =
(\zeta -z)\, {\rm tr}(\omega) ={\rm tr}( p_{n})\, {\rm tr}(\omega)$.
\end{proof}

We now proceed with the construction of our invariant. From the definition of $SB_n$, any element $\omega$ in $SB_n$ can be written as
$$
\omega_1^{\epsilon_1}\omega_2^{\epsilon_2}\ldots \omega_m^{\epsilon_m},
$$
 where $\omega_j\in\{\sigma_i, \tau_i\,;\,1\leq i\leq n-1\}$ and $\epsilon_i = +1$ or $-1$. If $\omega_j = \tau_j$ we set $\epsilon_j := +1$.

\begin{defn}\rm
 The {\it exponent} $\epsilon(\omega)$ of $\omega$ is defined as the sum
 $\epsilon_1 +\ldots +\epsilon_m$. Since $SB_n$ embeds in a group \cite{fkr}, $\epsilon(\omega)$ is well--defined.
\end{defn}

  Let now ${\mathcal X}_d =\{x_1,\ldots , x_{d-1}\}$ be a set satisfying the $E$--condition and let ${\mathcal S}$ be the set of oriented singular links. We define the following map on the set ${\mathcal S}$.

\begin{defn}\rm
Let $\omega\in SB_n$. We define the map $\Delta$ on the closure $\widehat{\omega}$ of $\omega$ as follows:
$$
\Delta (\widehat{\omega}) := \left(\frac{1-\lambda u}{\sqrt{\lambda}(1-u)\zeta}\right)^{n-1}
\left({\rm tr}\circ \delta_{\sqrt{\lambda}}\right)(\omega)
$$
where:
$$
 \lambda:= \frac{z -(1-u)\zeta}{uz}
$$
Equivalently, setting
$$
D:= \frac{1-\lambda u}{\sqrt{\lambda}(1-u)\zeta}
$$
 we can write:
$$
\Delta (\widehat{\omega}) = D^{n-1}(\sqrt{\lambda})^{\epsilon(\omega)}{\rm tr}(
\delta(\omega))
$$
For the definitions of $\delta_{\sqrt{\lambda}}$ and $\delta$ recall (\ref{deltaa}).
\end{defn}

\begin{thm}\label{invariant}
Assuming the $E$--condition, $\Delta$ is an isotopy invariant for oriented singular links.
\end{thm}

\begin{proof}
We need to show that $\Delta$ is well--defined on isotopy classes of oriented singular links. According to Theorem~\ref{gema}, it suffices to prove that $\Delta$ is consistent with moves (i), (ii) and (iii). From the facts that $\epsilon(\omega\omega^\prime) = \epsilon(\omega^\prime\omega)$ and
 ${\rm tr}(ab)= {\rm tr}(ba)$, it follows that $\Delta$ respects moves (i) and (ii).
Let now $\omega \in SB_n$. Then $\omega\sigma_n\in SB_{n+1}$ and $\epsilon(\omega\sigma_n)= \epsilon(\omega)+1$. Hence:
$$
\Delta(\widehat{\omega\sigma_n})
 = D^{n}(\sqrt{\lambda})^{\epsilon(\omega\sigma_n)}\, {\rm tr}(\delta(\omega \sigma_n))
 = D^{n}(\sqrt{\lambda})^{\epsilon(\omega)+1}\, {\rm tr}(\delta(\omega)g_n)
 = D \sqrt{\lambda}\, z \, \Delta(\widehat{\omega})
$$
where we used that ${\rm tr}(\delta (\omega)g_n) = z\, {\rm tr}(\delta (\omega))$.  Now $z = \frac{(1-u)\zeta}{1-\lambda u}$, so $D \sqrt{\lambda}\, z  = 1$. Therefore, $\Delta(\widehat{\omega\sigma_n}) = \Delta(\widehat{\omega}).$ Finally, we will prove that $\Delta(\widehat{\omega\sigma_n^{-1}}) = \Delta(\widehat{\omega})$.  Indeed:
$$
\Delta(\widehat{\omega\sigma_n^{-1}})
= D^{n}(\sqrt{\lambda})^{\epsilon(\omega\sigma_n^{-1})}{\rm tr}(\delta(\omega \sigma_n^{-1}))
= D^{n}(\sqrt{\lambda})^{\epsilon(\omega)-1}\, {\rm tr}(\delta(\omega)g_n^{-1})
$$
Resolving $g_n^{-1}$ from Eq.~\ref{invrs} we obtain:
$$
\Delta (\widehat{\omega\sigma_n^{-1}})
=  D^{n}(\sqrt{\lambda})^{\epsilon(\omega)-1}\,
\left[ z - (u^{-1}-1)\zeta + (u^{-1}-1)z \right]\, {\rm tr}(\delta(\omega))
$$
Also, from Theorem \ref{thmcon} and Eq. \ref{egomega} we have:
$$
{\rm tr}( \delta (\omega) e_{n}) = \zeta\,{\rm tr}(\delta (\omega))\quad \text{and} \quad
{\rm tr}( \delta (\omega) e_{n}g_n) = z\, {\rm tr}(\delta (\omega)).
$$
Therefore:
$$
\Delta (\widehat{\omega\sigma_n^{-1}})
=  D^{n}(\sqrt{\lambda})^{\epsilon(\omega)-1}\, \frac{z+(u-1)\zeta}{u}\, {\rm tr}(\delta(\omega))
= \frac{D}{\sqrt{\lambda}} \, \frac{z+(u-1)\zeta}{u}\,\Delta(\widehat{\omega})
= \Delta(\widehat{\omega})
$$
 Move (iii) of Theorem~\ref{gema} is now checked and the proof is concluded.
\end{proof}

\begin{rem}\label{remvas}
\rm
The invariant $\Delta$ is not of finite type. Indeed, take for example the link  $\widehat{\tau_1^k}$, which contains $k$ singular crossings. By the second relation of  Proposition~\ref{relp}, we have:
$\Delta(\widehat{\tau_1^k}) = (\sqrt{\lambda})^k (u+1)^{k-1}(\zeta- z)$, which is not equal to zero for all $k$. Of course, it would be interesting to consider an exponential variable change, and see if the coefficients become invariants of finite type as, for example, in the case of the Jones polynomial.
\end{rem}

\subsection {\it Skein relations}

Let $L_{+}$, $L_{-}$ and $L_{\times}$ be diagrams of three oriented singular links, which are identical, except near one crossing, where they are as follows:

\smallbreak

Then we have the following result.
\begin{prop}
The invariant $\Delta$ satisfies the following skein relation:
\begin{equation}\label{skein}
\frac{1}{\sqrt{\lambda}}\,\Delta (L_{+}) - \sqrt{\lambda}\,\Delta (L_{-})
= \frac{u^ {-1}-1}{\sqrt{\lambda}}\,\Delta (L_{\times})
\end{equation}
\end{prop}
\begin{proof}
The proof is standard. By the Alexander theorem for singular braids  we may assume that $L_+$ is in braided form and that $L_+ = \widehat{\beta \sigma_i}$ for some $\beta\in SB_n$.  Also that $L_- = \widehat{\beta \sigma_i^{-1}}$ and $L_{\times} = \widehat{\beta \tau_i}$. From the definition of $\Delta$ and by Theorem~\ref{invariant} we have:
$$
\frac{1}{\sqrt{\lambda}}\,\Delta (L_{+}) - \sqrt{\lambda}\,\Delta (L_{-})
=
D^{n-1}
(\sqrt{\lambda})^{\epsilon(\beta)}\left({\rm tr}(\delta(\beta\sigma_i)) - {\rm tr}(\delta(\beta\sigma_i^{-1}))\right)
$$
Now,
\begin{eqnarray*}
{\rm tr}(\delta(\beta\sigma_i)) - {\rm tr}(\delta(\beta\sigma_i^{-1}))
& = &  {\rm tr}(\delta(\beta)g_i) - {\rm tr}(\delta_1(\beta ) g_i^{-1})\\
& = &  {\rm tr}(\delta(\beta)(g_i -  g_i^{-1}))\\
& =  &  (u^{-1}-1)\,{\rm tr}(\delta(\beta)p_i) \qquad(\text{from\, Eq. \ref{gipi}})\\
 & = & (u^{-1}-1)\,{\rm tr}(\delta(\beta\tau_i))
 \end{eqnarray*}
Finally, substituting $(\sqrt{\lambda})^{\epsilon(\beta)}= (\sqrt{\lambda})^{-1} (\sqrt{\lambda})^{\epsilon(\beta\tau_i)}$ we deduce:
$$
\frac{1}{\sqrt{\lambda}}\,\Delta (L_{+}) - \sqrt{\lambda}\,\Delta (L_{-})
 = \frac{D^{n-1}}{\sqrt{\lambda}}
 (\sqrt{\lambda})^{\epsilon(\beta\tau_i)}(u^{-1}-1){\rm tr}(\delta(\beta)p_i)\\
 = \frac{u^ {-1}-1}{\sqrt{\lambda}}\,\Delta (L_{\times})
$$
Thus the proof is concluded.
\end{proof}

\subsection{\it Computations}

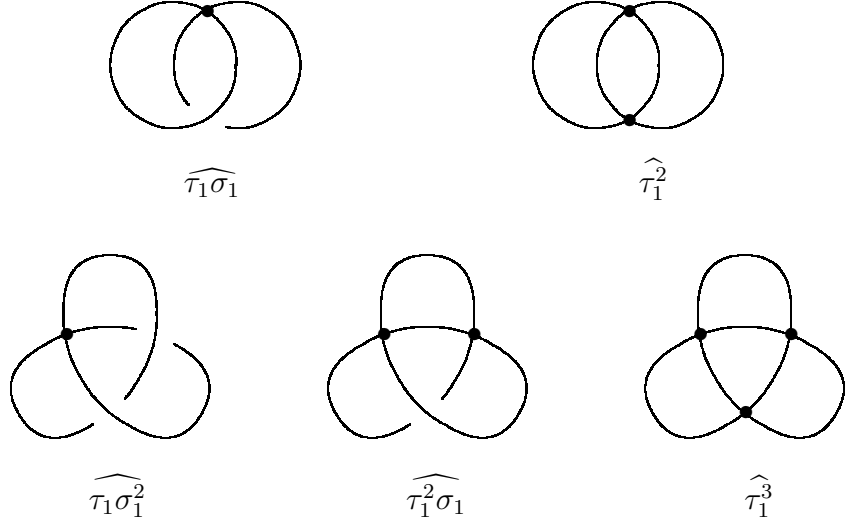
\begin{figure}[H]
\begin{center}
\setlength{\unitlength}{.8pt}
\begin{picture}(400,250) 

\put(92.5, 232){$\bullet$}
\put(292, 232){$\bullet$}
\put(292, 180){$\bullet$}
\qbezier(65.29,235.838)(74.443,241.514)(84.986,239.311)
\qbezier(52.123,220.338)(55.486,230.569)(64.979,235.658)
\qbezier(52,200)(48,209.999)(52.001,220)
\qbezier(64.979,184.342)(55.487,189.43)(52.124,199.663)
\qbezier(84.986,180.69)(74.444,178.486)(65.29,184.163)
\qbezier(102.66,190.751)(96,182.288)(85.341,180.751) 
\qbezier(109.31,209.82)(110.07,199.056)(102.892,191.027)
\qbezier(102.891,228.973)(110.07,220.945)(109.732,210.18)
\qbezier(85.34,239.248)(96,237.713)(102.66,229.249) 
\qbezier(124.71,235.838)(115.557,241.514)(105.014,239.311)
\qbezier(137.877,220.338)(134.514,230.569)(125.021,235.658)
\qbezier(138,200)(142,209.999)(137.999,220)
\qbezier(125.021,184.342)(134.513,189.43)(137.876,199.663)
\qbezier(105.014,180.69)(115.556,178.486)(124.71,184.163)
\qbezier(80.269,209.82)(79.93,199.056)(87.108,191.027)
\qbezier(87.109,228.973)(79.93,220.945)(80.268,210.18)
\qbezier(104.66,239.248)(94,237.713)(87.34,229.249) 


\qbezier(265.29,235.838)(274.443,241.514)(284.986,239.311)
\qbezier(252.123,220.338)(255.486,230.569)(264.979,235.658)
\qbezier(252,200)(248,209.999)(252.001,220)
\qbezier(264.979,184.342)(255.487,189.43)(252.124,199.663)
\qbezier(284.986,180.69)(274.444,178.486)(265.29,184.163)
\qbezier(302.66,190.751)(296,182.288)(285.341,180.751) 
\qbezier(309.31,209.82)(310.07,199.056)(302.892,191.027)
\qbezier(302.891,228.973)(310.07,220.945)(309.732,210.18)
\qbezier(285.34,239.248)(296,237.713)(302.66,229.249) 
\qbezier(324.71,235.838)(315.557,241.514)(305.014,239.311)
\qbezier(337.877,220.338)(334.514,230.569)(325.021,235.658)
\qbezier(338,200)(342,209.999)(337.999,220)
\qbezier(325.021,184.342)(334.513,189.43)(337.876,199.663)
\qbezier(305.014,180.69)(315.556,178.486)(324.71,184.163)
\qbezier(287.34,190.751)(294,182.288)(304.659,180.751) 
\qbezier(280.269,209.82)(279.93,199.056)(287.108,191.027)
\qbezier(287.109,228.973)(279.93,220.945)(280.268,210.18)
\qbezier(304.66,239.248)(294,237.713)(287.34,229.249) 


\put(85,150){$\widehat{\tau_1\sigma_1}$}
\put(300,150){$\widehat{\tau_1^2}$}


\put(26,79){$\bullet$}
\put(176,79){$\bullet$}
\put(218.5,79){$\bullet$}
\put(325.5,79){$\bullet$}
\put(368.5,79){$\bullet$}
\put(347,42){$\bullet$}

\qbezier(50,120)(73,120)(72,92)
\qbezier(50,86)(58,86)(62,85)
\qbezier(50,120)(27,120)(28,92)
\qbezier(50,86)(42,86)(38,85)
\qbezier(6.699,45)(-4.801,64.919)(19.947,78.053)
\qbezier(36.144,62)(32.144,68.928)(31.01,72.892)
\qbezier(6.699,45)(18.199,25.081)(41.947,39.947)
\qbezier(36.144,62)(40.144,55.072)(43.01,52.108)
\qbezier(93.301,45)(81.801,25.081)(58.053,39.947)
\qbezier(63.856,62)(59.856,55.072)(56.99,52.108)
\qbezier(93.301,45)(104.801,64.919)(80.053,78.053)
\qbezier(63.856,62)(67.856,68.928)(68.99,72.892)
\qbezier(68.99,72.892)(72.517,83)(72,92.001) 
\qbezier(38,85)(27.483,83)(19.947,78.052) 
\qbezier(31.01,72.892)(27.483,83)(28,92.001) 

\qbezier(43.01,52.108)(50,44)(58.053,39.947) 


\qbezier(200,120)(223,120)(222,92)
\qbezier(200,86)(208,86)(212,85)
\qbezier(200,120)(177,120)(178,92)
\qbezier(200,86)(192,86)(188,85)
\qbezier(156.699,45)(145.199,64.919)(169.947,78.053)
\qbezier(186.144,62)(182.144,68.928)(181.01,72.892)
\qbezier(156.699,45)(168.199,25.081)(191.947,39.947)
\qbezier(186.144,62)(190.144,55.072)(193.01,52.108)
\qbezier(243.301,45)(231.801,25.081)(208.053,39.947)
\qbezier(213.856,62)(209.856,55.072)(206.99,52.108)
\qbezier(243.301,45)(254.801,64.919)(230.053,78.053)
\qbezier(213.856,62)(217.856,68.928)(218.99,72.892)
\qbezier(218.99,72.892)(222.517,83)(222,92.001) 
\qbezier(212,85)(222.517,83)(230.053,78.052) 
\qbezier(188,85)(177.483,83)(169.947,78.052) 
\qbezier(181.01,72.892)(177.483,83)(178,92.001) 
\qbezier(193.01,52.108)(200,44)(208.053,39.947) 


\qbezier(350,120)(373,120)(372,92)
\qbezier(350,86)(358,86)(362,85)
\qbezier(350,120)(327,120)(328,92)
\qbezier(350,86)(342,86)(338,85)
\qbezier(306.699,45)(295.199,64.919)(319.947,78.053)
\qbezier(336.144,62)(332.144,68.928)(331.01,72.892)
\qbezier(306.699,45)(318.199,25.081)(341.947,39.947)
\qbezier(336.144,62)(340.144,55.072)(343.01,52.108)
\qbezier(393.301,45)(381.801,25.081)(358.053,39.947)
\qbezier(363.856,62)(359.856,55.072)(356.99,52.108)
\qbezier(393.301,45)(404.801,64.919)(380.053,78.053)
\qbezier(363.856,62)(367.856,68.928)(368.99,72.892)
\qbezier(368.99,72.892)(372.517,83)(372,92.001) 
\qbezier(362,85)(372.517,83)(380.053,78.052) 
\qbezier(338,85)(327.483,83)(319.947,78.052) 
\qbezier(331.01,72.892)(327.483,83)(328,92.001) 
\qbezier(343.01,52.108)(350,44)(358.053,39.947) 
\qbezier(356.99,52.108)(350,44)(341.947,39.947) 

\put(40,0){$\widehat{\tau_1\sigma_1^2}$}
\put(190,0){$\widehat{\tau_1^2\sigma_1}$}
\put(350,0){$\widehat{\tau_1^3}$}

\end{picture}
\caption{Examples of singular knots and links}\label{fig5}
\end{center}
\end{figure}

In this subsection we compute the values of the invariant $\Delta$ on some basic classical and singular knots, assuming always the $E$--condition. The singular ones are illustrated in Figure~\ref{fig5}. We shall first give some formulas that are useful for computations.
For powers of $g_i$ we can easily deduce by induction the following formulae.
\begin{lem}\label{powers}
Let $m \in {\Bbb Z}, k \in {\Bbb N}$. (i) For $m$ positive, define:

\noindent $\alpha_m = (u-1)
\sum_{l=0}^{k-1}u^{2l}$ if $m=2k$ \ and \ $\beta_m = u(u-1)
\sum_{l=0}^{k-1}u^{2l}$ if $m=2k+1$. Then:
\[
g_i^{m} =\left\{
\begin{array}{l} 1 + \alpha_m p_i \quad \text{ \ if \ } m=2k \\
 g_i - \beta_mp_i \quad \text{ \ if \ } m=2k+1
\end{array}\right.
\]
(ii) For $m$ negative, define:

\noindent $\alpha'_m = u^{-1}(u^{-1}-1) \sum_{l=0}^{k-1}u^{-2l}$ if $m=-2k$
 \ and \ $\beta'_m = (u^{-1}-1) \sum_{l=0}^{k-1}u^{-2l}$ if $m=-2k+1$.
Then:
\[
g_i^{m} = \left\{
\begin{array}{l} 1 + \alpha'_mp_i \quad \text{ \ if \ } m=-2k \\
 g_i - \beta'_mp_i \quad \text{ \ if \ } m=-2k+1
\end{array}\right.
\]
\end{lem}

\smallbreak

We now proceed with our computations.
\smallbreak

\noindent $\bullet$ Clearly, for the unknot ${\rm O}$, $\Delta ({\rm O})=1$.

\smallbreak
\noindent $\bullet$  Let ${\rm K}_1=\widehat{\tau_1}$. Then $e(\tau_1)=1$, so
$
\Delta ({\rm K}_1) = D \sqrt{\lambda}\, {\rm tr}(p_1) =D \sqrt{\lambda}\,[{\rm tr}(e_1)- {\rm tr}(e_1g_1)]
= D \sqrt{\lambda}\,[\zeta- z]$. Then:
$$
\Delta({\rm K}_1) = \frac{\zeta- z}{z}
$$

\smallbreak
\noindent $\bullet$ Let ${\rm H}=\widehat{\sigma_1^2}$, the Hopf link. We have ${\rm tr }(g_1^2) = {\rm tr }(1 + (u+1)p_1)= 1 + (u+1)(\zeta - z)$ and $e(\sigma_1^2) = 2$. Then:
$$
\Delta({\rm H}) =\frac{1-\lambda u}{(1-u)\zeta}\sqrt{\lambda}\, \left(1 + (u+1)(\zeta - z)\right) = z^{-1}\sqrt{\lambda}
\left(1 + (u+1)(\zeta - z)\right).
$$

\smallbreak
 \noindent $\bullet$ Let ${\rm H}_1=\widehat{\sigma_1\tau_1}$. We have ${\rm tr}(g_1p_1)= -u{\rm tr}(p_1)= -u(\zeta - z)$ and $e(\sigma_1\tau_1) = 2$. Then:
$$
\Delta({\rm H}_1)= -\frac{1-\lambda u}{(1-u)\zeta}\sqrt{\lambda} \,u(\zeta - z) =
z^{-1}\sqrt{\lambda}u(\zeta - z).
$$

\smallbreak
\noindent $\bullet$ Let ${\rm H}_2=\widehat{\tau_1^2}$. We have $e(\tau_1^2) = 2$. Then $\Delta ({\rm H}_2) =D \lambda\, {\rm tr}(p_1^2) = D \lambda\,(u+1) {\rm tr}(p_1)$. So:
$$
\Delta({\rm H}_2) =\frac{1-\lambda u}{(1-u)\zeta}\sqrt{\lambda}
(u + 1)(\zeta- z) = z^{-1}\sqrt{\lambda}(u + 1)(\zeta- z).
$$

\smallbreak
\noindent $\bullet$ Let ${\rm T}=\widehat{\sigma_1^3}$, the right-handed trefoil.  We have  $g_1^3 = g_1 -u(u-1)e_1 + u(u-1)e_1g_1$ from (Lemma~\ref{powers}). Hence:
${\rm tr}(g_1^3) = z -u(u-1)\zeta + u(u-1)z$. Moreover $e(\sigma_1^ 3)= 3$.
Then, using that $1-\lambda u = z^ {-1}\zeta(1-u)$, we obtain:
\begin{eqnarray*}
\Delta({\rm T}) & = & D(\sqrt{\lambda})^3\left[(u(u-1) + 1)z-u(u-1)\zeta\right]\\
& = &
\frac{\lambda}{z}\left[(u(u-1) + 1)z-u(u-1)\zeta\right].
\end{eqnarray*}

\smallbreak
\noindent $\bullet$
Let $T_1=\widehat{\tau_1\sigma_1^2}$. We have $p_1g_1^2= -up_1g_1= u^ 2p_1$ so, ${\rm tr}(p_1g_1^2) = u^2(\zeta -z)$. Moreover, $e(\tau_1\sigma_1^2)= 3$. Then
$$
\Delta ({\rm T}_1) = D (\sqrt{\lambda})^3{\rm tr}(p_1g_1^ 2) = \frac{u^2\lambda}{z}(\zeta -z).
$$

\smallbreak
\noindent $\bullet$ Let ${\rm T}_2=\widehat{\tau_1^ 2\sigma_1}$. We have $p_1^ 2g_1= -u(u + 1)p_1$ so, ${\rm tr}(p_1^2g_1) = -u(u+1)(\zeta -z)$. Moreover, $e(\tau_1^ 2\sigma_1)= 3$. Then:
$$
\Delta ({\rm T}_2) = D (\sqrt{\lambda})^3{\rm tr}(p_1^2g_1) = \frac{-u(u+1)\lambda}{z}(\zeta -z).
$$

\smallbreak
\noindent $\bullet$ Let ${\rm T}_3=\widehat{\tau_1^3}$. We have $p_1^3= (u + 1)^2p_1$, so ${\rm tr}(p_1^3) = (u+1)^2(\zeta -z)$. Moreover, $e(\tau_1^3)= 3$. Then:
$$
\Delta ({\rm T}_3) = D (\sqrt{\lambda})^3{\rm tr}(p_1^3) = \frac{(u+1)^2\lambda}{z}(\zeta -z).
$$


\end{document}